\documentclass[
    ,final            
,numberedheadings 
  ]
  {aipproc}

\layoutstyle{8x11double}

\usepackage{amsfonts,amssymb,amsmath}

 \newtheorem{thm}{Theorem}[section]
 \newtheorem{lem}[thm]{Lemma}
 \newtheorem{defn}[thm]{Definition}
 
\newcommand{\be}{\begin{equation}}
\newcommand{\ee}{\end{equation}}

\newcommand{\cl}[1]{\ensuremath{Cl(#1)}} 

\newcommand{\vect}[1]{\ensuremath{\mbox{\textbf{\textit{#1}}}}}

  \newcommand{\field}[1]{\mathbb{#1}}
  \newcommand{\C}{\field{C}}
  \newcommand{\R}{\field{R}}  
  \newcommand{\Z}{\field{Z}} 

  \newcommand{\N}{\field{N}}
  
  \newcommand{\HQ}{\field{H}}
  \renewcommand*{\H}{\ensuremath{\mathbb{H}}}
  \newcommand{\qi}{\ensuremath{\mbox{\boldmath $i$}}}
  \newcommand{\sqi}{\ensuremath{\mbox{\boldmath $\scriptstyle i$}}}
  \newcommand{\qj}{\ensuremath{\mbox{\boldmath $j$}}}
  \newcommand{\sqj}{\ensuremath{\mbox{\boldmath $\scriptstyle j$}}}
  \newcommand{\qk}{\ensuremath{\mbox{\boldmath $k$}}}

\newcommand{\M}{\mathcal{M}}

\newcommand{\sandwich}{(\,)} 

\begin{document}

\title{Clifford Fourier-Mellin transform with two real square roots of $-1$ in $Cl(p,q)$, $p+q=2$}

\classification{02.30.Nw,07.05.Pj,02.10.Hh}

\keywords{Fourier-Mellin transform, Clifford algebra, square roots of $-1$ }

\author{Eckhard Hitzer}{
  address={University of Fukui, 910-8507 Fukui, Japan}
}

\begin{abstract}
We describe a non-commutative generalization of the complex Fourier-Mellin transform to Clifford algebra valued signal functions over the domain $\R^{p,q}$ taking values in $Cl(p,q)$, $p+q=2$.
\end{abstract}

\maketitle


\section{Introduction}

The Fourier-Mellin transform is an excellent tool in order to achieve translation, rotation and scale invariant shape recognition. It has recently been extended to the algebra of quaternions $\mathbb{H}$ in the form of the quaternionic Fourier-Mellin Transform \cite{EH:QFMT}, which in principle also allows a non-marginal processing of color images. Beyond this new progress has been made in the study and complete classification of square roots of $-1$ in real Clifford algebras $Cl(p,q)$ \cite{SJS:Biqroots,HA:GeoRoots-1,HHA:ICCA9,worksheets}. In this paper\footnote{Please respect the Creative Peace License regarding the content of this paper \cite{EH:CPL}.} we make therefore a first attempt to generalize the Fourier-Mellin transform to Clifford algebras $Cl(p,q)$, $p+q=2$, using a set of two real square roots of $-1$, $f,g\in Cl(p,q), f^2=g^2=-1$. This includes the case of quaternions, because of the isomorphism $\mathbb{H} \equiv Cl(0,2)$, but in addition the algebras $Cl(2,0)$ and $Cl(1,1)$ are also included. First, for the quaternion Fourier transform a \textit{split} of quaternions with respect to any two pure unit quaternions has been developed \cite{EH:QFTgen,EH:OPS-QFT}. It has also been applied to the quaternionic Fourier-Mellin Transform, and it has been generalized to Clifford algebras \cite{EH:AGACSE2012}. We will also apply this split.

\section{Clifford's geometric algebra}
\label{sc:CliffAlg}

\begin{defn}[Clifford's geometric algebra \cite{FM:ICNAAM2007,PL:CAandSpin}] \label{df:CliffAlg}
Let  $\{  e_1, e_2, \ldots , e_p, e_{p+1}, \ldots$, $e_n  \}$, with $n=p+q$, $e_k^2=\varepsilon_k$, $\varepsilon_k = +1$ for $k=1, \ldots , p$, $\varepsilon_k = -1$ for $k=p+1, \ldots , n$, be an \textit{orthonormal base} of the inner product vector space $\R^{p,q}$ with a geometric product according to the multiplication rules 
\be
  e_k e_l + e_l e_k = 2 \varepsilon_k \delta_{k,l}, 
  \qquad k,l = 1, \ldots n,
\label{eq:mrules}
\ee
where $\delta_{k,l}$ is the Kronecker symbol with $\delta_{k,l}= 1$ for $k=l$, and $\delta_{k,l}= 0$ for $k\neq l$. This non-commutative product and the additional axiom of \textit{associativity} generate the $2^n$-dimensional Clifford geometric algebra $Cl(p,q) = Cl(\R^{p,q}) = Cl_{p,q} = \mathcal{G}_{p,q} = \R_{p,q}$ over $\R$. The set $\{ e_A: A\subseteq \{1, \ldots ,n\}\}$ with $e_A = e_{h_1}e_{h_2}\ldots e_{h_k}$, $1 \leq h_1< \ldots < h_k \leq n$, $e_{\emptyset}=1$, forms a graded (blade) basis of $Cl(p,q)$. The grades $k$ range from $0$ for scalars, $1$ for vectors, $2$ for bivectors, $s$ for $s$-vectors, up to $n$ for pseudoscalars. 
The vector space $\R^{p,q}$ is included in $Cl(p,q)$ as the subset of 1-vectors. The general elements of $Cl(p,q)$ are real linear combinations of basis blades $e_A$, called Clifford numbers, multivectors or hypercomplex numbers.
\end{defn}

In general $\langle A \rangle_{k}$ denotes the grade $k$ part of $A\in Cl(p,q)$. The parts of grade $0$ and $k+s$, respectively, of the geometric product of a $k$-vector $A_k\in Cl(p,q)$ with an $s$-vector $B_s\in Cl(p,q)$ 
\begin{gather}
  A_k \ast B_s := \langle A_k B_s \rangle_{0}, 
  \qquad
  A_k \wedge B_s := \langle A_k B_s \rangle_{k+s},
  \label{eq:gaprods}
\end{gather}
are called \textit{scalar product} and \textit{outer product}, respectively.

For Euclidean vector spaces $(n=p)$ we use $\R^{n}=\R^{n,0}$ and $Cl(n) = Cl(n,0)$. Every $k$-vector $B$ that can be written as the outer product $B = \vect{b}_1 \wedge \vect{b}_2 \wedge \ldots \wedge \vect{b}_k$ of $k$ vectors $\vect{b}_1, \vect{b}_2, \ldots, \vect{b}_k \in \R^{p,q}$ is called a \textit{simple} $k$-vector or \textit{blade}. 

Multivectors
  $M \in \cl{p,q}$ have $k$-vector parts ($0\leq k \leq n$):
  {scalar} part
  $Sc(M) = \langle M \rangle = \langle M \rangle_0 = M_0 \in \R$, 
  {vector} part
  $\langle M \rangle_1 \in \R^{p,q}$, 
  {bi-vector} part
  $\langle M \rangle_2$,  \ldots, 
  and
  {pseudoscalar} part $\langle M \rangle_n\in\bigwedge^n\R^{p,q}$
\begin{equation}\label{eq:MVgrades}
    M  =  \sum_{A} M_{A} \vect{e}_{A}
       =  \langle M \rangle + \langle M \rangle_1 + \langle M \rangle_2 + \ldots +\langle M \rangle_n \, .
\end{equation}

The \textit{principal reverse} of $M \in \cl{p,q}$ defined as
\begin{equation}\label{eq:MVrev}
  \widetilde{M}=\; \sum_{k=0}^{n}(-1)^{\frac{k(k-1)}{2}}\langle \overline{M} \rangle_k,
\end{equation}
often replaces {complex conjugation and quaternion conjugation}. Taking the \textit{reverse} is equivalent to reversing the order of products of basis vectors in the basis blades $e_A$. The operation $\overline{M}$ means to change in the basis decomposition of $M$ the sign of every vector of negative square $\overline{e_A} = \varepsilon_{h_1}e_{h_1}\varepsilon_{h_2}e_{h_2}\ldots \varepsilon_{h_k}e_{h_k}$, $1 \leq h_1< \ldots < h_k \leq n$. Reversion, $\overline{M}$, and principal reversion are all involutions. 

The principal reverse of every basis element $e_A \in \cl{p,q}$, $1 \leq A \leq 2^n$, has the property
\be 
  \widetilde{e_A} \ast e_B = \delta_{AB},\qquad 1 \leq A,B \leq 2^n ,
\ee 
where the Kronecker delta $\delta_{AB}= 1$ if $A=B$, and $\delta_{AB}= 0$ if $A\neq B$. For the vector space $\R^{p,q}$ this leads to a reciprocal basis $e^l$, $1 \leq l,k \leq n$
\begin{gather}
  e^l:= \widetilde{e_l}=\varepsilon_l e_l,
  \nonumber \\
  e^l \ast e_{k} = e^l \cdot e_k 
  =\left\{ 
  \begin{array}{ll}
  1,   & \quad \text{for}\,\,\, l = k \\
  0,   & \quad \text{for}\,\,\, l \neq k 
  \end{array}
  \right.
  .
\end{gather}

For $M,N \in \cl{p,q}$ we get $M\ast \widetilde{N}=\sum_{A} M_A N_A.$
  Two multivectors $M,N \in \cl{p,q}$ are \textit{orthogonal} if and only if $M\ast \widetilde{N} = 0$.
  The {modulus} $|M|$ of a multivector $M \in \cl{p,q}$ is defined as 
  \be
     |M|^2 = {M\ast\widetilde{M}}= {\sum_{A} M_A^2}.
  \ee

\subsection{Multivector signal functions}

A multivector valued function
  $f: \R^{p,q} \rightarrow \cl{p,q}$, has $2^n$ blade components
  $(f_A: \R^{p,q} \rightarrow \R)$
  \begin{equation}\label{eq:MVfunc}
    f(\mbox{\textbf{\textit{x}}})  =  \sum_{A} f_{A}(\vect{x}) {\vect{e}}_{A},
    \qquad
    \vect{x} = \sum_{l=1}^{n} x_l e^l = \sum_{l=1}^{n} x^l e_l. 
  \end{equation}
We define the \textit{inner product} of two
 functions  $f, g : \R^{p,q} \rightarrow \cl{p,q}$ by
\begin{align}
  \label{eq:mc2}
  (f,g) 
  &= \int_{\R^{p,q}}f(\vect{x})
    \widetilde{g(\vect{x})}\;d^n\vect{x}
  \nonumber \\
  &= \sum_{A,B}\vect{e}_A \widetilde{\vect{e}_B}
    \int_{\R^{p,q}}f_A (\vect{x})
    g_B (\vect{x})\;d^n\vect{x},
\end{align}
with the \textit{symmetric scalar part}
\begin{align}
  \label{eq:symsc}
  \langle f,g\rangle 
  &= \int_{\R^{p,q}}f(\vect{x})\ast\widetilde{g(\vect{x})}\;d^n\vect{x}
  \nonumber \\
  &= \sum_{A}
    \int_{\R^{p,q}}f_A (\vect{x})g_A (\vect{x})\;d^n\vect{x},
\end{align}
and the $L^2(\mathbb{R}^{p,q};\cl{p,q})$-\textit{norm} 
\begin{align}
   \label{eq:0mc2}
   \|f\|^2 
   = \left\langle ( f,f ) \right\rangle
   &= \int_{\R^{p,q}} |f(\vect{x})|^2 d^n\vect{x} 
   \nonumber \\
   &= \sum_{A} \int_{\R^{p,q}} f_A^2(\vect{x})\;d^n\vect{x},
\end{align}
\begin{align}
   L^2(&\R^{p,q};\cl{p,q}) 
   \nonumber \\
   &= \{f: \R^{p,q} \rightarrow \cl{p,q} \mid \|f\| < \infty \}. 
\end{align}

The \textit{vector derivative} $\nabla$ of a function $f : \R^{p,q} \rightarrow \cl{p,q}$ can be expanded in a basis of $\R^{p,q}$ as \cite{GS:ConfMpsGA} 
\be 
  \nabla = \sum_{l=1}^n e^l \partial_l
  \quad \text{with} \quad \partial_l=\partial_{x_l}=\frac{\partial}{\partial x_l}, 
  \,\,\, 1 \leq l \leq n .
\ee

\subsection{Square roots of $-1$ in Clifford algebras}

\begin{figure}
    \includegraphics[height=.3\textheight]{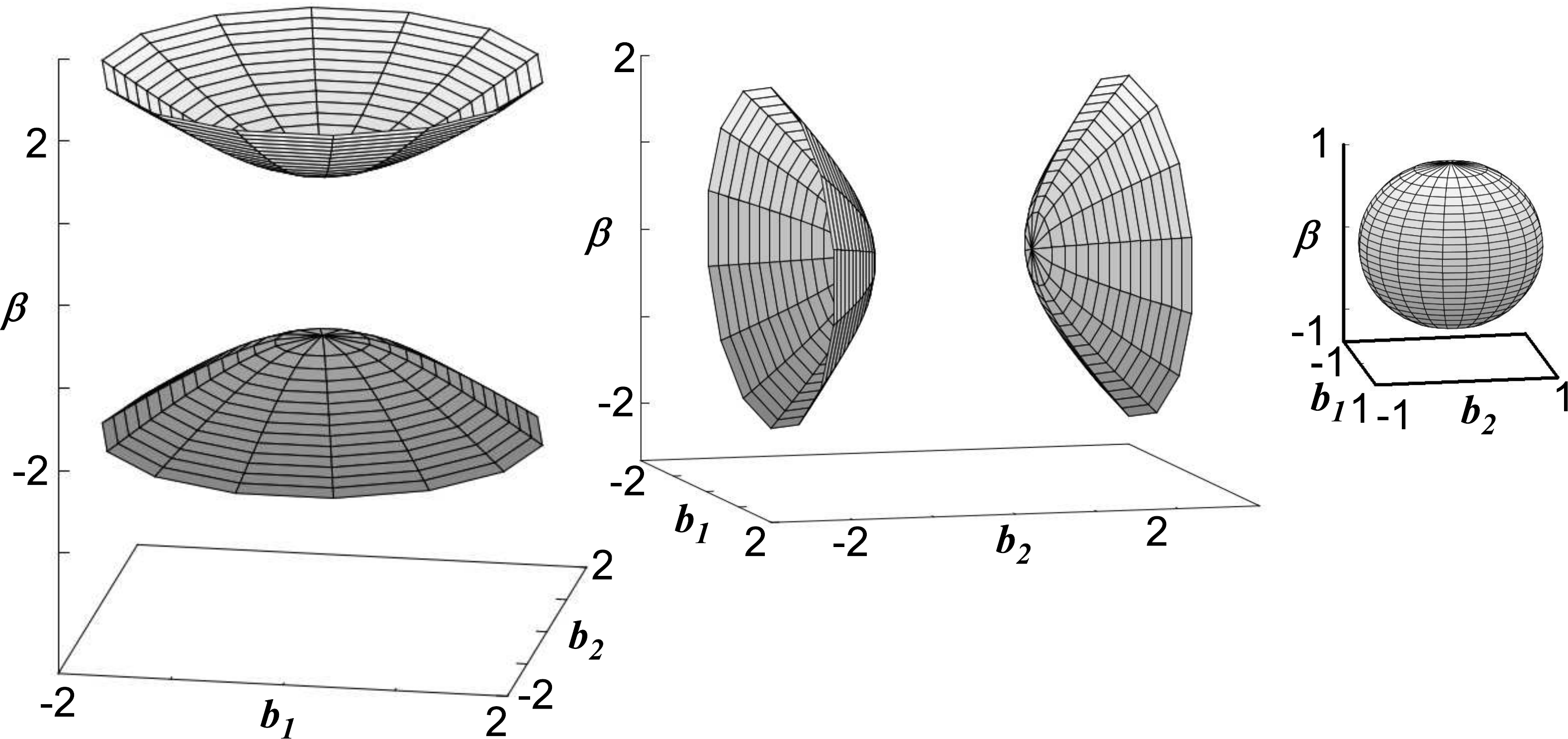}
    \caption{Manifolds \cite{HHA:ICCA9} of square roots $f$ of $-1$ in $\cl{2,0}$ (left), $\cl{1,1}$ (center), and $\cl{0,2}\cong\H$ (right). The square roots are $f=\alpha + b_1 e_1+b_2e_2+\beta e_{12},$ with $\alpha, b_1, b_2, \beta \in \R$, $\alpha=0$, and $\beta^2=b_1^2e_2^2+b_2^2e_1^2+e_1^2e_2^2$.}  
    \label{fg:Cln=2}
\end{figure}

\subsubsection{The $\pm$ split with respect to two square\\ roots of $-1$}
\label{sc:pmsplit}

With respect to any square root $f\in \cl{p,q}$ of $-1$, $f^2=-1$, every multivector $A\in \cl{p,q}$ can be split into \textit{commuting} and \textit{anticommuting} parts~\cite{HHA:ICCA9}.
\begin{lem}
\label{lm:fsplit}
Every multivector $A \in \cl{p,q}$ has, with respect to a square root $f\in \cl{p,q}$ of~$-1$, i.e., $f^{-1}=-f,$ the unique decomposition
\begin{gather} 
  A_{+f} = \frac{1}{2}(A + f^{-1}Af),\qquad  
  A_{-f} = \frac{1}{2}(A - f^{-1}Af)
  \nonumber \\
  A = A_{+f}+A_{-f}, 
  \qquad A_{+f}\,f = f A_{+f}, 
  \nonumber \\
  A_{-f}\,f = -fA_{-f}. 
\end{gather}
\end{lem}

For $f,g \in \cl{p,q}$ an arbitrary pair of square roots of $-1$, $f^2=g^2=-1$, the map $f\sandwich g$ is an involution, because $f^2 x g^2 = (-1)^2x = x, \forall \, x\in \cl{p,q}$. In \cite{EH:QFTgen} a split of quaternions by means of the pure unit quaternion basis elements $\vect{i}, \vect{j} \in \HQ$ was introduced, and generalized to a general pair of pure unit quaternions in \cite{EH:OPS-QFT}. We here use a \textit{generalized} version of the split for $\cl{p,q}$ \cite{EH:AGACSE2012}.

\begin{defn}[$\pm$ split with respect to two $\sqrt{-1}$]
\label{df:genOPS}
Let $f,g$ $\in$ $\cl{p,q}$ be an arbitrary pair of square roots of $-1$, $f^2=g^2=-1$, including the cases $f = \pm g$. The general $\pm$ split is then defined with respect to the two square roots $f, g$ of $-1$ as 
\begin{equation}
  \label{eq:genOPS}
  x_{\pm} = \frac{1}{2}(x \pm f x g), \qquad \forall\, x\in \cl{p,q}.
\end{equation} 
\end{defn}

Note that the split of Lemma \ref{lm:fsplit} is a special case of Definition \ref{df:genOPS} with $g=-f$. 

We observe from (\ref{eq:genOPS}), that $f x g = x_+ - x_-$, i.e. under the map $f\sandwich g$ the $x_+$ part is invariant, but the $x_-$ part changes sign
\begin{equation}
  fx_{\pm}g 
  = \pm x_{\pm}.
  \label{eq:fgrotqm}
\end{equation} 

The two parts $x_{\pm}$ can be represented with Lemma \ref{lm:fsplit} as linear combinations of $x_{+f}$ and $x_{-f}$, or of $x_{+g}$ and $x_{-g}$
\begin{align} 
  x_{\pm} 
  &= x_{+f}\,\frac{1\pm fg}{2} + x_{-f}\,\frac{1\mp fg}{2}
  \nonumber \\
  &= \frac{1\pm fg}{2}\,x_{+g} + \frac{1\mp fg}{2}\,x_{-g} .
\end{align} 

For $Cl(0,2)\cong \M(2d,\C)$, or for both $f,g$ being blades in $Cl(2,0)$ or $Cl(1,1)$, we have $\widetilde{f}=-f$, $\widetilde{g}=-g$. We therefore obtain the following lemma. 

\begin{lem}[Orthogonality of two $\pm$ split parts \cite{EH:AGACSE2012}]
  \label{lm:OPSortho}
  Given any two multivectors $x,y \in \cl{p,q}$ and applying the $\pm$ split (\ref{eq:genOPS}) with respect to two square roots of $-1$ we get zero for the
  scalar part of the mixed products
  \begin{equation}
     Sc(x_+\widetilde{y_-}) = 0, 
     \qquad Sc(x_-\widetilde{y_+}) = 0 .
  \end{equation}
\end{lem}

Finally, we have the \textit{general identity} \cite{EH:AGACSE2012}
\be 
  e^{\alpha f} x_{\pm} e^{\beta g} 
  = x_{\pm} e^{(\beta\mp\alpha) g}
  = e^{(\alpha \mp \beta)f } x_{\pm} .
  \label{eq:expqexp}
\ee

\section{The Clifford Fourier Mellin transformations (CFMT)}

\subsection{Robert Hjalmar Mellin (1854--1933)}

\begin{figure}
  \includegraphics[height=3.0cm]{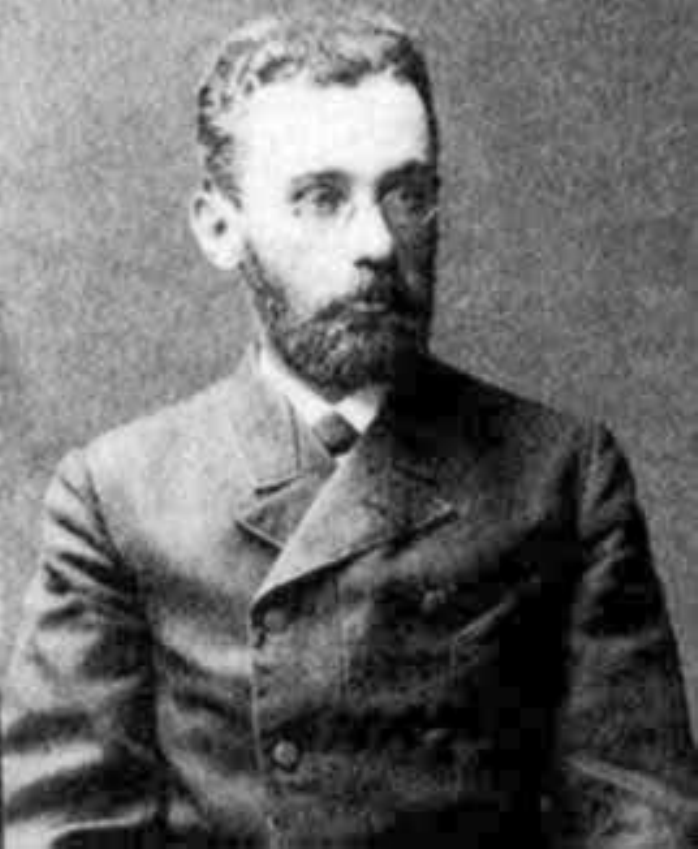}
  \caption{Robert Hjalmar Mellin (1854--1933). Image: \cite{StA:Mellinbio} and Wikipedia. 
  \label{fg:RHMellin}}
\end{figure}

Robert Hjalmar Mellin (1854--1933) \cite{StA:Mellinbio}, Fig. \ref{fg:RHMellin}, was a Finnish mathematician, a student of G. Mittag-Leffler and K. Weierstrass. He became the director of the Polytechnic Institute in Helsinki, and in 1908 first professor of mathematics at Technical University of Finland. He was a fervent fennoman with fiery temperament, and co-founder of the Finnish Academy of Sciences. He became known for the \textit{Mellin transform} with major applications to the evaluation of integrals, see \cite{PBM:EvIntMellT}, which lists {1624 references}. During his last 10 years he tried to {refute Einstein's theory} of relativity as logically untenable. 

\begin{defn}[Classical Fourier Mellin transf. (FMT)]
\begin{align}
  &\forall(v,k)\in \mathbb{R} \times \mathbb{Z}, \quad
  \mathcal{M}\{h\}(v,k)= 
  \nonumber \\
  &= \frac{1}{2\pi}\int_0^{\infty}\int_0^{2\pi}
  h(r,\theta) r^{-iv}e^{-ik\theta}d\theta \frac{dr}{r},
\end{align}
where $h: \mathbb{R}^2 \rightarrow \R$ denotes a function representing, e.g., a gray level image defined over a compact set of $\mathbb{R}^2$. 
\end{defn}
Well known applications are to shape recognition (independent of rotation and scale), image registration, and similarity. 

\subsection{Inner product, symmetric part, norm of Clifford-valued functions}

For $Cl(0,2)\cong \M(2d,\C)$, or for both $f,g$ being blades in $Cl(2,0)$ or $Cl(1,1)$, we have $\widetilde{f}=-f$, $\widetilde{g}=-g$. We therefore obtain the following Pythagorean modulus identity for the $L^2(\R^2; Cl(p,q))$-norm, $p+q=2$, i.e., 
\be
  \|h\|^2 = \|h_+\|^2 + \|h_-\|^2 \,\,.
\ee

We now define the generalization of the FMT to Clifford-valued signals. 

\begin{defn}[Clifford FMT (CFMT)]
Let $f,g \in Cl(p,q): f^2=g^2=-1$, $p+q=2$, be any pair of real square roots of $-1$ in $Cl(p,q)$, $p+q=2$. 
The Clifford Fourier Mellin transform (CFMT) is given by
\begin{align}
  &\forall(v,k)\in \mathbb{R} \times \mathbb{Z},\quad
  \hat{h}(v,k)=
  \mathcal{M}\{h\}(v,k) =
  \nonumber \\
  &= \frac{1}{2\pi}\int_0^{\infty}\int_0^{2\pi}
  r^{-f v} h(r,\theta) e^{-g k\theta}d\theta \frac{dr}{r},
\end{align}
where $h: \mathbb{R}^2 \rightarrow Cl(p,q)$ denotes a function from $\mathbb{R}^2$ into the real Clifford $Cl(p,q)$, $p+q=2$, 
such that $|h|$ is summable over $\mathbb{R}_+^* \times \mathbb{S}^1$ under the measure $d\theta \frac{dr}{r}$. $\mathbb{R}_+^*$ is the multiplicative group of positive and non-zero real numbers.
\end{defn}
For special pure unit quaternion (isomorphic to $Cl(0,2)$) values $f=\qi$, $g=\qj$ we have the special case
\begin{align}
  \forall(v,k)\in \mathbb{Z}\times \mathbb{R},\quad
  \hat{h}(v,k)=
  \mathcal{M}\{h\}(v,k) =
  \nonumber \\  
  = \frac{1}{2\pi}\int_0^{\infty}\int_0^{2\pi}
  r^{-\sqi v} h(r,\theta) e^{-\sqj k\theta}d\theta \frac{dr}{r},
\end{align}
Note, that the $\pm$ split and the CFMT {commute}: 
$$
  \mathcal{M}\{h_{\pm}\} = \mathcal{M}\{h\}_{\pm}.
$$

\begin{thm}[Inverse CFMT]
The CFMT can be inverted by
\begin{align} 
  h(r,\theta) &= 
  \mathcal{M}^{-1}\{h\}(r,\theta)
  \nonumber \\
  &= \frac{1}{2\pi}\int_{-\infty}^{\infty}\sum_{k\in \Z}
  r^{f v} \,\hat{h}(v,k) \,e^{g k\theta} dv.
\end{align}
\end{thm}
The proof uses
\begin{gather}
  \frac{1}{2\pi}\sum_{k\in \Z}e^{g k(\theta-\theta')} 
  = \delta(\theta-\theta'),
  \qquad
  r^{f v} = e^{f v \ln r}, 
  \nonumber \\
  \frac{1}{2\pi}\int_0^{2\pi}e^{f v (\ln (r) - s)} dv 
  = \delta(\ln (r) - s) .
\end{gather}

We now investigate the basic properties of the CFMT. 
First, left linearity:
For $\alpha, \beta \in \{q \mid q = q_r + q_f f, \, q_r, q_f \in \R\}$,
\begin{gather} 
  m(r,\theta) = \alpha h_1(r,\theta) + \beta h_2(r,\theta)
  \nonumber \\
  \Longrightarrow \,\,\, 
  \hat{m}(v,k) = \alpha \hat{h}_1(v,k) + \beta \hat{h}_2(v,k) . 
\end{gather}

Second, right linearity:
For $\alpha', \beta' \in \{q \mid q = q_r + q_g g, \, q_r, q_g \in \R\}$,
\begin{gather}
  m(r,\theta) = h_1(r,\theta)\alpha' + h_2(r,\theta)\beta'
  \nonumber \\
  \Longrightarrow \,\,\, 
  \hat{m}(v,k) = \hat{h}_1(v,k)\alpha' + \hat{h}_2(v,k)\beta' . 
\end{gather}

The linearity of the CFMT leads to
\begin{align}
  \mathcal{M}\{h\}(v,k)
  &= \mathcal{M}\{h_- + h_+\}(v,k)
  \nonumber \\
  &= \mathcal{M}\{h_-\}(v,k)
    + \mathcal{M}\{ h_+\}(v,k),
\end{align}
which gives rise to the following thoerem about \textit{the quasi-complex FMT like forms for CFMT of} $h_{\pm}$.

\begin{thm}[Quasi-complex forms for CFMT]
\label{th:fpmtrafo}
The CFMT of $h_{\pm}$ parts of $h \in L^2(\R^2,\HQ)$ 
have simple {quasi-complex forms}
\begin{align}
 \mathcal{M}\{h_{\pm}\} 
  &{=} \frac{1}{2\pi}\int_0^{\infty}\int_0^{2\pi}
    h_{\pm} r^{\pm g v} e^{-g k \theta} d\theta \frac{dr}{r}
  \nonumber \\
  &{=} \frac{1}{2\pi}\int_0^{\infty}\int_0^{2\pi}
    r^{-f v} e^{\pm f k \theta} h_{\pm} d\theta \frac{dr}{r} \,\, .
\end{align}
\end{thm}
Theorem \ref{th:fpmtrafo} allows to use \textit{discrete} and \textit{fast} software to compute the CFMT based on a pair of complex FMT transformations. 

For $Cl(0,2)\cong \M(2d,\C)$, or for both $f,g$ being blades in $Cl(2,0)$ or $Cl(1,1)$, we have $\widetilde{f}=-f$, $\widetilde{g}=-g$, and under these conditions we have for the two split parts of the CFMT, the following lemma. 
\begin{lem}[Modulus identities]
Due to $|x|^2 = |x_-|^2 + |x_+|^2$, for $Cl(0,2)\cong \M(2d,\C)$, or for both $f,g$ being blades in $Cl(2,0)$ or $Cl(1,1)$, we get for $f:\R^2\rightarrow Cl(p,q)$, $p+q=2$, the following identities
\begin{align}
  |h(r,\theta)|^2 
  &= |h_-(r,\theta)|^2 + |h_+(r,\theta)|^2,
  \nonumber \\
  |\mathcal{M}\{h\}(v,k)|^2 
  &= |\mathcal{M}\{h_-\}(v,k)|^2 + |\mathcal{M}\{h_+\}(v,k)|^2.
\end{align}
\end{lem}

Further properties are \textit{scaling} and \textit{rotation}:
For 
$ m(r,\theta) = h(ar,\theta+\phi)$, 
$\,a > 0, \,\, 0\leq \phi \leq 2\pi $,
\be
  \widehat{m}(v,k) = a^{f v} \hat{h}(v,k) e^{g k \phi} .
\ee 

Moreover, we have the following magnitude identity:
\be 
  |\widehat{m}(v,k)| = |\hat{h}(v,k)| ,
  \label{eq:magid}
\ee 
i.e. the magnitude of the CFMT of a scaled and rotated quaternion signal 
$m(r,\theta) = h(ar,\theta+\phi)$ is identical to the magnitude of the CFMT of $h$.
Equation (\ref{eq:magid}) forms the basis for applications to {rotation and scale invariant shape recognition} and image registration. This may now be extended to {color images}, since quaternions can encode colors RGB in their $\qi, \qj, \qk$ components, and to signals with values in $Cl(2,0)$ and $Cl(1,1)$.

The \textit{reflection at the unit circle} ($r \rightarrow \frac{1}{r}$) leads to
\begin{gather}
  m(r,\theta) = h(\frac{1}{r}, \theta) 
  \nonumber \\
  \Longrightarrow \qquad
  \widehat{m}(v,k) = \hat{h}(-v,k) .
\end{gather}

\textit{Reversing} the sense of sense of \textit{rotation} ($\theta \rightarrow -\theta$) yields
\begin{gather}
  m(r,\theta) = h(r, -\theta) 
  \nonumber \\
  \Longrightarrow \qquad
  \widehat{m}(v,k) = \hat{h}(v,-k) .
\end{gather}

Regarding radial and rotary modulation we assume
\be
  m(r,\theta) = r^{f v_0} \,h(r,\theta) \,e^{g k_0\theta}, 
  \qquad
  v_0 \in \R, \, k_0 \in \Z . 
\ee
Then we get
\be
  \widehat{m}(v,k) = \hat{h}(v-v_0, k-k_0) .
\ee

\subsection{CFMT derivatives and power scaling}

We note for the logarithmic derivative that $\frac{d}{d \ln r} = r \frac{d}{dr} = r \partial_r$, 
\be 
  \mathcal{M}\{(r\partial_r)^n h\}(v,k) = (f v)^n \hat{h}(v,k), 
  \qquad n \in \N .
\ee
Applying the angular derivative with respect to $\theta$ we obtain
\be 
  \mathcal{M}\{\partial_{\theta}^n h\}(v,k) = \hat{h}(v,k) (g k)^n , 
  \qquad n \in \N .
\ee
Finally, power scaling with $\ln r$ and $\theta$ leads for all $m,n \in \N$, to
\be 
  \mathcal{M}\{ (\ln r)^m \theta^n h\}(v,k)
  = f^m \,\partial_v^m \partial_k^n \hat{h}(v,k) \,g^n 
  .
\ee

\subsection{CFMT Plancherel and Parseval theorems}

For the CFMT we have the following two theorems. 

\begin{thm}[CFMT Plancherel theorem]
The scalar part of the inner product of two functions $h,m : \R^2 \rightarrow Cl(p,q)$, $p+q=2$,
is
\be
  \langle h,m \rangle = \langle \hat{h},\widehat{m} \rangle .
\ee 
\end{thm}
\begin{thm}[CFMT Parseval theorem]
Let $h : \R^2 \rightarrow Cl(p,q)$, $p+q=2$, and assume $Cl(0,2)\cong \M(2d,\C)$, or for both $f,g$ being blades in $Cl(2,0)$ or $Cl(1,1)$, then
\be
  \|h\| = \|\hat{h}\| ,
  \qquad
  \|h\|^2 = \|\hat{h}\|^2 = \|\hat{h}_+\|^2 + \|\hat{h}_-\|^2 .
\ee 
\end{thm}

\section{Symmetry of the CFMT}

The CFMT of real signals analyzes symmetry. 
The following notation will be used\footnote{In this section we assume $g\neq \pm f$, but a similar study is possible for $g=\pm f$.}. 
The function $h_{ee}$ is \textit{even} with respect to (w.r.t.)
$r \rightarrow \frac{1}{r} 
\Longleftrightarrow 
\ln r \rightarrow -\ln r$,
i.e. w.r.t. the reflection at the unit circle,
and \textit{even} w.r.t. $\theta \rightarrow -\theta$,
i.e. w.r.t. reversing the sense of rotation
(reflection at the $\theta = 0$ line of polar coordinates
in the ($r,\theta$)-plane). Similarly we denote by
$h_{eo}$ even-odd symmetry, by 
$h_{oe}$ odd-even symmetry, and by 
$h_{oo}$ odd-odd symmetry.   

Let $h$ be a real valued function $\R^2 \rightarrow \R$. The CFMT of $h$ results in
\be
  \hat{h}(v,k) 
  = \underbrace{\hat{h}_{ee}(v,k)}_{\text{real part}}
    + \underbrace{\hat{h}_{eo}(v,k)}_{f\text{-part}}
    + \underbrace{\hat{h}_{oe}(v,k)}_{g\text{-part}}
    + \underbrace{\hat{h}_{oo}(v,k)}_{fg\text{-part}} \, .
\ee 
The CFMT of a real signal therefore automatically separates components
with different combinations of symmetry w.r.t. reflection at the unit circle and
reversal of the sense of rotation.
The four components of the CFMT kernel differ by radial and angular phase shifts.

\section{Conclusions}

We have generalized the Fourier-Mellin transform to a Clifford Fourier Mellin transform (CFMT) acting on $Cl(p,q)$-valued signals, with $p+q=2$, which includes the previously treated case of quaternions \cite{EH:QFMT}, as well as the algebras $Cl(2,0)$ and $Cl(1,1)$. We have derived several properties of this transform: inversion, linearity, quasi-complex (split) forms, and a (split) modulus identity. Beyond this for scaling, rotation, the total modulus remains invariant. This forms the foundation for translation, rotation and scale invariant shape recognition. We have further studied symmetry properties of the CFMT, derivatives and power scaling, Plancherel and Parseval theorems.


\begin{theacknowledgments}
  \textit{Soli Deo Gloria.} E.H. thanks his family for their patient support of his work. He further thanks S. J. Sangwine, R. Ablamowicz, J. Helmstetter, J. Morao, S. Svetlin, und W. Sproessig, for good cooperation, support and advice. 
\end{theacknowledgments}



\bibliographystyle{aipproc}   


\end{document}